\documentclass{amsart}

\usepackage{latexsym}

\usepackage{amsfonts,amssymb,euscript,epsfig,color}

\usepackage[all]{xy}

\usepackage[english,francais]{babel}

\newtheorem{theorem}{Theorem}[section]
\newtheorem{lemma}[theorem]{Lemma}
\newtheorem{proposition}[theorem]{Proposition}

\newcommand{\hgot}{\ensuremath{\mathfrak{h}}}
\newcommand{\kgot}{\ensuremath{\mathfrak{k}}}

\newcommand{\tgot}{\ensuremath{\mathfrak{t}}}

\newcommand{\Ccal}{\ensuremath{\mathcal{C}}}

\newcommand{\Fcal}{\ensuremath{\mathcal{F}}}

\newcommand{\Lcal}{\ensuremath{\mathcal{L}}}

\newcommand{\Ocal}{\ensuremath{\mathcal{O}}}

\newcommand{\Scal}{\ensuremath{\mathcal{S}}}

\newcommand{\Z}{\ensuremath{\mathbb{Z}}}
\newcommand{\C}{\ensuremath{\mathbb{C}}}
\newcommand{\Lbb}{\ensuremath{\mathbb{L}}}

\newcommand{\R}{\ensuremath{\mathbb{R}}}

\newcommand{\QS}{\ensuremath{\mathrm{Q}^{\mathrm{spin}}}}

\newcommand{\mm}{\ensuremath{\hbox{\rm m}}}

\begin{document}

\title{The multiplicities of the equivariant index of twisted   Dirac operators}

\author{Paul-Emile  PARADAN, Mich\`ele VERGNE}

\address{Institut de Math\'ematiques et de Mod\'elisation de Montpellier (I3M), UMR CNRS 5149, Universit\'e Montpellier 2} 
\email{Paul-Emile.Paradan@math.univ-montp2.fr}

\address{Institut de Math\'ematiques de Jussieu, UMR CNRS 7586, Universit\'e Paris-Diderot paris 7} 
\email{Vergne@math.jussieu.fr}

\maketitle

\begin{abstract}

In this note, we give a geometric expression for the multiplicities of the equivariant index of a Dirac operator twisted by a line bundle.

\end{abstract}

\vskip 0.5\baselineskip

\section{Introduction}

This note is an announcement of work whose details will appear later.

Let $M$ be a compact connected manifold. We assume that $M$ is even dimensional and oriented.
We consider a spin$^c$ structure on $M$, and denote by $\mathcal S$ the corresponding irreducible Clifford module.
Let $K$ be a compact connected Lie group acting on $M$, and preserving the spin$^c$ structure.
We denote by $D: \Gamma(M, \Scal^+)\to \Gamma(M, \Scal^-)$ the corresponding twisted Dirac operator.
The equivariant index  of $D$, denoted $\QS_K(M)$,  belongs to the Grothendieck group of representations of $K$,
$$
\QS_K(M) = \sum_{\pi\in \widehat{K}} \mm(\pi) \ \pi.
$$

An important example is when $M$ is a compact complex manifold, $K$ a compact group of holomorphic transformations of $M$, and $\Lcal$
any holomorphic $K$-equivariant line bundle on $M$ (not necessarily ample). Then  the Dolbeaut operator twisted by $\Lcal$ can be
realized as a twisted Dirac operator $D$. In this case  $\QS_K(M)= \sum_q (-1)^q H^{0,q}(M,\Lcal)$.

The aim of this note is to give a geometric description of the multiplicity $\mm(\pi)$ in the spirit of  the  Guillemin-Sternberg phenomenon $[Q,R]=0$ \cite{Guillemin-Sternberg82,Meinrenken98,meinrenken-sjamaar, Tian-Zhang98, pep-RR}.

Consider the determinant line bundle $\Lbb= \det(\Scal)$ of the  spin$^c$ structure. This is a $K$-equivariant complex line bundle on $M$.
The choice of a $K$-invariant hermitian metric  and of a $K$-invariant hermitian connection  $\nabla$ on $\Lbb$ determines an abstract moment map
$$\Phi_{\nabla} : M\to \kgot^*$$
by the relation $\Lcal(X)-\nabla_{X_M}= \frac{i}{2}\langle \Phi_{\nabla},X\rangle$, for  all   $X\in\kgot.$
We  compute $\mm(\pi)$ in term of the reduced ``manifolds" $\Phi_{\nabla}^{-1}(f)/K_f$.
  This formula extends the result of \cite{pep-spin}. However, in this note, we do not assume any hypothesis on the line bundle $\Lbb$, in particular we do not assume that the curvature of the connection $\nabla$ is a symplectic form.
In this pre-symplectic setting, a (partial) answer to this question has been obtained by \cite{Karshon-Tolman93,Grossberg-Karshon94,Grossberg-Karshon98,CdS-K-T} when $K$ is a torus.
Our method is based on localization techniques as in \cite{pep-RR}, \cite{pep-spin}.

\section{Admissible coadjoints orbits}

We consider  a compact connected Lie group  $K$ with Lie algebra $\kgot$. Consider an admissible  coadjoint orbit $\Ocal$ (as in \cite{du}), oriented by its symplectic structure.
Then $\Ocal$  carries a
$K$-equivariant bundle of spinors $\Scal_\Ocal$, such that the associated moment
map  is the injection $\Ocal$ in $ \kgot^*$. We denote
by $\QS_K(\Ocal)$ the corresponding equivariant index.

Let us  describe the admissible coadjoint orbits with their spin$^c$ index.

Let $T$ be a Cartan subgroup of $K$ with Lie algebra $\mathfrak t$.
Let $\Lambda\subset {\mathfrak t}^*$  be the lattice  of weights of $T$ (thus $e^{i\lambda}$ is a character of $T$).
Choose a positive system $\Delta^+\subset {\mathfrak t}^*$, and let
$\rho=\frac{1}{2}\sum_{\alpha\in \Delta^+}\alpha$.
Let $ {\mathfrak t}^*_{\geq 0}$ be the closed Weyl chamber and  we denote by $\mathcal{F}$  the set of the relative interiors of the  faces of ${\mathfrak t}^*_{\geq 0}$.
 Thus ${\mathfrak t}^*_{\geq 0}=\coprod_{\sigma\in {\mathcal F}} \sigma$, and we denote $ {\mathfrak t}^*_{>0}\in \mathcal{F}$ the interior of ${\mathfrak t}^*_{\geq 0}$.

We index the set $\hat K$ of  classes of finite dimensional irreducible representations of $K$ by
the set $(\Lambda+\rho)\cap {\mathfrak t}^*_{> 0}$.
The irreducible representation $\pi_\lambda$ corresponding to
$\lambda\in (\Lambda+\rho)\cap {\mathfrak t}^*_{>0}$ is the irreducible representation with
infinitesimal character $\lambda$. Its highest weight is $\lambda-\rho$.

Let  $\sigma\in \mathcal{F}$. The stabilizer $K_\xi$ of a point $\xi\in \sigma$ depends only of $\sigma$.
We denote it by $K_\sigma$, and by $\kgot_\sigma$ its Lie algebra. We choose on ${\mathfrak k}_\sigma$
the system  of positive roots contained in  $\Delta^+$, and let $\rho_\sigma$ be the corresponding $\rho$.

When $\mu\in\sigma$, the coadjoint orbit $K\cdot \mu$
is admissible  if and only if
 $\mu-\rho+\rho_\sigma \in \Lambda$.
The spin$^c$ equivariant index of the admissible orbits is  described in the following lemma.

\begin{lemma}Let $K\cdot\mu$ be an admissible orbit: $\mu\in\sigma$ and $\mu-\rho+\rho_\sigma \in \Lambda$.  If
$\mu+\rho_\sigma$ is regular, then $\mu+\rho_\sigma\in \rho+\overline{\sigma}$. Thus  we have
$$
\QS_K(K\cdot\mu)=
\begin{cases}
   0\qquad\qquad {\rm if}\ \mu+\rho_\sigma\ \mathrm{is\ singular},\\
   \pi_{\mu+\rho_\sigma}\qquad \, {\rm if}\ \mu+\rho_\sigma\ \mathrm{is\ regular}.
\end{cases}
$$
In particular, if
$\lambda\in (\Lambda+\rho)\cap {\mathfrak t}^*_{>0}$, then
$K\cdot\lambda$ is admissible and
$\QS_K(K\cdot\lambda)=\pi_\lambda$.

\end{lemma}

Let $\mathcal{H}_{\mathfrak k}$ be the set of conjugacy classes of the reductive algebras ${\kgot}_f,f\in\mathfrak k^*$. We denote by
$\mathcal S_{\mathfrak k}$ the set of conjugacy classes of the semi-simple parts $[\hgot,\hgot]$ of the elements $(\hgot)\in \mathcal H_{\mathfrak k}$.
The map $(\hgot)\to( [\hgot,\hgot])$ induces a bijection between $\mathcal H_{\mathfrak k}$ and $\mathcal S_{\mathfrak k}$.

The map $\Fcal\longrightarrow \mathcal{H}_\kgot$, $\sigma\mapsto (\kgot_\sigma)$, is surjective and for $({\mathfrak h})\in\mathcal H_{\mathfrak k}$ we denote by

$\bullet$  ${\mathcal F}({\mathfrak h})$ the set of $\sigma\in {\mathcal F}$ such that $(\mathfrak k_\sigma)= (\mathfrak h)$,

$\bullet$ $\mathfrak k^*_{\mathfrak h}\subset \mathfrak k^*$ the set of elements $f\in\kgot^*$   with infinitesimal stabilizer $\kgot_f$  belonging to the conjugacy class $({\mathfrak h})$.

We have  $\mathfrak k^*_{\mathfrak h}=K\left(\cup_{\sigma \in {\mathcal F}(\mathfrak h)}\sigma\right)$. In particular all coadjoint orbits contained in  $\mathfrak k^*_{\mathfrak h}$ have the same dimension.
We say that such a coadjoint orbit is of type $(\hgot)$. If $(\hgot)=(\tgot)$, then  $\mathfrak k^*_{\mathfrak h}$ is the open subset of regular elements.

We  denote by $A(\hgot)$ the set of admissible coadjoint orbits  of type $(\hgot)$. This is a discrete subset of orbits in
 $\mathfrak k^*_{\mathfrak h}$.

{\bf Example 1 :}  Consider the group $K=SU(3)$ and let $(\hgot)$ be the conjugacy class such that $\kgot^*_\hgot$
is equal to the set of subregular element $f\in \kgot^*$ (the orbit of $f$ is of dimension $\dim (K/T)-2$).
Let $\omega_1,\omega_2$ be the two fundamental weights.
Let $\sigma_1,\sigma_2$ be the half lines $\R_{>0}\omega_1$, $\R_{>0}\omega_2$.
Then $\kgot^*_{\hgot}\cap \tgot^*_{\geq 0}=\sigma_1\cup \sigma_2$.
The set $A(\hgot)$ is equal to the collection of orbits $K\cdot (\frac{1+2n}{2}\omega_i), n\in\Z_{\geq 0}, i=1,2$. The
representation $\QS_K(K\cdot (\frac{1+2n}{2}\omega_i))$ is $0$ is $n=0$, otherwise it is
the irreducible representation $\pi_{\rho+(n-1)\omega_i}$.
In particular, both representations associated to the admissible orbits $\frac{3}{2}\omega_1$ and $\frac{3}{2}\omega_2$
are the trivial representation $\pi_\rho$.

\section{The theorem}

Consider the action of $K$ in $M$. Let $(\kgot_M)$
be the conjugacy class of the generic infinitesimal stabilizer. On a $K$-invariant open and dense subset of $M$, the conjugacy class of
$\kgot_m$ is  equal to  $(\kgot_M)$.
Consider the (conjugacy class)  $([\kgot_M,\kgot_M])$.

We start by stating two vanishing lemmas.
\begin{lemma}
If $([\kgot_M,\kgot_M])$ does not belong to the set $\Scal_\kgot$, then  $\QS_K(M)=0$ for any $K$-invariant spin$^c$ structure on $M$.
\end{lemma}

If $([\kgot_M,\kgot_M])=([{\mathfrak h},{\mathfrak h}])$ for some ${\mathfrak h}\in \mathcal H_{\kgot}$, any $K$-invariant map $\Phi: M\to {\mathfrak k}^*$ is such that $\Phi(M)$  is included in  the closure of ${\mathfrak k}^*_{\mathfrak h}$.
\begin{lemma}
Assume that $([{\mathfrak k}_M,{\mathfrak k}_M])=([{\mathfrak h},{\mathfrak h}])$ with $\hgot\in \mathcal H_\mathfrak k$.
Let us consider a spin$^c$ structure on $M$ with determinant bundle $\Lbb$.
If there exists a $K$-invariant hermitian connection $\nabla$ on $\Lbb$ such that
$\Phi_{\nabla}(M)\cap {\mathfrak k}^*_\mathfrak h=\emptyset$, then
$\QS_K(M)=0$.
\end{lemma}

Thus from now on, we assume that  the action of $K$ on $M$ is such that
$([\kgot_M,\kgot_M])=([{\mathfrak h},{\mathfrak h}])$ for some ${\mathfrak h}\in \mathcal H_{\kgot}$.
Let us consider a spin$^c$ structure on $M$ with determinant bundle $\Lbb$ and a $K$-invariant hermitian connection with moment map $\Phi_{\nabla}:M\to \kgot^*.$

We extend the definition of the index to disconnected even dimensional oriented manifolds by defining $\QS_K(M)$ to be the sum over the connected components of $M$.
If $K$ is the trivial group, $\QS_K(M)\in \Z$ and is denoted simply by $\QS(M)$.

Consider a coadjoint orbit $\Ocal=K\cdot f$. The reduced space
$M_\Ocal$ is defined to be the topological space $\Phi_{\nabla}^{-1}(\Ocal)/K=\Phi_{\nabla}^{-1}(f)/K_f$.
We also denote it by $M_f$.
This space might not be connected.

In the next section, we define a $\Z$-valued function $\Ocal\mapsto\QS(M_\Ocal)$ on the set $A(\hgot)$ of admissible
orbits of type $(\hgot)$. We call it the reduced index :

$\bullet$ if $M_\Ocal=\emptyset$,  then $\QS(M_\Ocal)=0$,

$\bullet$ when $M_\Ocal$ is an orbifold, the reduced index $\QS(M_\Ocal)$ is defined as an index of a Dirac operator associated to a natural ``reduced" spin$^c$ structure on $M_\Ocal$.

Otherwise, it is defined via a limit procedure. Postponing this definition, we have the following theorem.

\begin{theorem}\label{QR}
Assume that $([{\mathfrak k}_M,{\mathfrak k}_M])=([{\mathfrak h},{\mathfrak h}])$ with $(\hgot)\in \mathcal H_\mathfrak k$.
Then
$$
\QS_K(M)=\sum_{\Ocal\in A(\hgot)}  \QS(M_{\Ocal}) \ \QS_K(\Ocal).
$$
\end{theorem}

In the expression above, when $\hgot$ is not abelian,
$\QS_K(\Ocal)$ can be $0$, and  several orbits
 $\Ocal\in A(\hgot)$ can give the same representation.

Theorem \ref{QR} is in the spirit  of the $[Q,R]=0$ theorem.
However it has some radically new features.
First, as $\Phi_{\nabla}$ is not the moment map of a Hamiltonian structure, the definition of the reduced space requires more care. For example, the fibers of $\Phi_{\nabla}$ might not be connected, and
 the Kirwan set  $\Phi_{\nabla}(M)\cap \tgot^*_{\geq 0}$ is not a convex polytope. Furthermore, this Kirwan set  depends  of the choice of connection $\nabla$.
Second, the map  $\Ocal\in A(\hgot)\to \QS_K(\Ocal)$ is not injective, when  $\hgot$ is not abelian. Thus the multiplicities
$\mm_\lambda$  of the representation $\pi_\lambda$ in
$\QS_K(M)$
will be eventually  obtained as  a sum of reduced indices involving several reduced spaces.

We explicit this last point.

\begin{theorem}\label{th:mult}
Assume that $([{\mathfrak k}_M,{\mathfrak k}_M])=([{\mathfrak h},{\mathfrak h}])$ with $(\hgot)\in \mathcal H_\mathfrak k$.
Let  $\mm_\lambda\in\Z$ be the multiplicity of  the  representation $\pi_\lambda$ in $\QS_K(M)$.
We have
\begin{equation}\label{eq:m-lambda}
\mm_\lambda=\sum_{\stackrel{\sigma\in \Fcal(\hgot)}{\lambda-\rho_\sigma\in \sigma}} \QS(M_{\lambda-\rho_\sigma}).
\end{equation}
\end{theorem}

More explicitly, the sum is taken over the  (relative interiors of) faces $\sigma$ of the Weyl chamber such that
\begin{equation}\label{eq:condition}
([\kgot_M,\kgot_M])=([\kgot_\sigma,\kgot_\sigma]),\hspace{0.5cm}\Phi_{\nabla}(M)\cap \sigma\neq \emptyset,\hspace{0.5cm}\lambda\in\{\sigma +\rho_\sigma\}.
\end{equation}

If  $\kgot_M$ is abelian, we have simply
 $\mm_\lambda= \QS( \Phi_{\nabla}^{-1}(\lambda)/T)$.
In particular, if the group $K$ is the circle group, and $\lambda$ is a regular value of the moment map $\Phi_{\nabla}$, this result was obtained in \cite{CdS-K-T}.

If $\kgot_M$ is not abelian,
 and the curvature of the connection $\nabla$ is symplectic, Kirwan convexity theorem implies that the image $\Phi_{\nabla}(M)\cap \tgot^*_{\geq 0}$ is contained in the closure of one single $\sigma$. Thus there is a unique $\sigma$ satisfying Conditions (\ref{eq:condition}). In this setting Theorem \ref{th:mult} is obtained in \cite{pep-spin}.

Let us give an example where several $\sigma$ contribute to the multiplicity of a representation $\pi_\lambda$.

We take the notations of Example 1. We label $\omega_1,\omega_2$ so that   $\kgot_{\omega_1}$ is the
group $S(U(2)\times U(1))$ stabilizing the line $\C e_3$  in the fundamental representation of $SU(3)$
in $\C^3=\C e_1\oplus \C e_2\oplus \C e_3$.

Let $P=\{0\subset L_2\subset L_3\subset \C^4\}$ be the partial flag manifold  with $L_2$ a subspace of $\C^4$ of dimension $2$ and $L_3$ a subspace of $\C^4$ of dimension $3$. Denote by $\Lcal_1,\Lcal_2$ the equivariant  line bundles on $P$ with fiber at $(L_2,L_3)$ the one-dimensional spaces  $\wedge^2L_2$ and $ L_3/L_2$ respectively.
Let $M$ be the subset of $P$ where $L_2$ is assumed to be a subspace of $\C^3$. 
Thus $M$ is fibered over $P_2(\C)$ with fiber $P_1(\C)$.
The group $SU(3)$ acts naturally on $M$, and the generic stabilizer of the action is $SU(2)$.
 We denote by   $\Lcal_{a,b}$ the line bundle
 $\Lcal_1^a \otimes \Lcal_2^{b}$ restricted to $M$. This line bundle is equipped with a natural  holomorphic and hermitian connection $\nabla$.
 Consider the spin$^c$ structure with  determinant  bundle $\Lbb=\Lcal_{2a+1,2b+1}$, where $a,b$ are positive integers.
If $a\geq b$, the curvature of the line bundle $\Lbb$ is non degenerate, and we are in the symplectic case.
Let us consider $b>a$.
It is easy to see that, in this case,  the Kirwan set $\Phi_\nabla(M)\cap \tgot^*_{\geq 0}$ is the non convex set  $[0,
b-a]\omega_1\cup [0,a+1]\omega_2$.
We compute the character of the representation  $\QS_K(M)$
by the Atiyah-Bott fixed point formula, and find
$$
\QS_K(M)=
\sum_{j=0}^{b-a-2} \pi_{\rho+j \omega_1} \oplus
\sum_{j=0}^{a-1} \pi_{\rho+j \omega_2}.
$$
In particular the multiplicity of $\pi_\rho$ (the trivial representation) is equal to $2$. We use now Theorem \ref{QR} and the discussion of Example 1, and   obtain  (reduced multiplicities are equal to $1$)
 $$
 \QS_K(M)=\sum_{j=0}^{b-a-1}\QS_K(K\cdot(\frac{1+2j}{2}  \omega_1))\oplus
 \sum_{j=0}^{a}  \QS_K(K\cdot(\frac{1+2j}{2} \omega_2)).
 $$
Using the  formulae  for $\QS_K(K\cdot(\frac{1+2n}{2} \omega_i))$ given in Example 1, these two formulae (fortunately) coincide.
Furthermore we see that both faces $\sigma_1,\sigma_2$  give a non zero contribution to  the multiplicity of the trivial representation.

\section{Definition of the reduced index}
We start by defining the reduced index for the action of an abelian torus  $H$ on a connected  manifold $Y$. Denote by $\Lambda$ the lattice of weights of $H$.
We do not assume $Y$ compact, but we assume that  the set  of stabilizers $H_m$  of points in $Y$ is finite.
Let $\hgot_Y$  be  the generic infinitesimal stabilizer of the action $H$ on $Y$, and $H_Y$ be  the  connected subgroup of $H$ with Lie algebra $\hgot_Y$. Thus $H_Y$ acts trivially on $Y$.
Let us consider a spin$^c$ structure on $Y$ with determinant bundle $\Lbb$, and a $H$ invariant connection $\nabla$ on $\Lbb$. The image $\Phi_\Delta(Y)$ spans an affine space $I_Y$ parallel to $\hgot_Y^{\perp}$. We assume that the fibers of the map $\Phi_\Delta$ are compact.
We can easily prove that there exists a finite collection of hyperplanes $W^1,\ldots, W^p$ in $I_Y$ such
that the group $H/H_Y$ acts locally freely on
$\Phi_{\Delta}^{-1}(f)$, when $f$ is in  $\Phi_\nabla(Y)$, but  not on any  of the hyperplanes $W^i$.

\begin{proposition}

$\bullet$ When $\mu\in I_Y\cap \Lambda$ is a regular value of $\Phi_\nabla: Y\to I_Y$, the reduced space
$Y_\mu$ is an oriented orbifold equipped with an induced spin$^c$ structure: we denote $\QS(Y_\mu)$ the corresponding spin$^c$ index.

$\bullet$ For any connected component $\Ccal$ of
$I_Y\setminus \cup_{k=1}^pW^k$,  we can associate a periodic polynomial function
$q^\Ccal : \Lambda\cap I_Y \to \Z$ such that
$$
q^\Ccal(\mu)= \QS(Y_\mu)
$$
for any element $\mu\in \Lambda\cap \Ccal$ which is a regular value of $\Phi: Y\to I_Y$.

$\bullet$ If $\mu\in\Lambda$ belongs to the closure of two connected components $\Ccal_1$ and $\Ccal_2$ of $I_Y\setminus \cup_{k=1}^pW^k$, we have
$$
q^{\Ccal_1}(\mu)=q^{\Ccal_2}(\mu).
$$
\end{proposition}


We can now state the definition of the ``reduced'' index  on $\Lambda$:

$\bullet$ $\QS(Y_\mu)=0$ if $\mu\notin \Lambda\cap I_Y$,

$\bullet$ for any $\mu\in \Lambda\cap I_Y$, we define $\QS(Y_\mu)$ as being equal to
$q^{\Ccal}(\mu)$ where $\Ccal$ is any connected component containing $\mu$ in its closure.
In fact $\QS(Y_\mu)$ is computed as an index of a particular spin$^c$ structure on the orbifold $\Phi_\nabla^{-1}(\mu+\epsilon)/H$ for any $\epsilon$ small and such that $\mu+\epsilon$ is a regular value of $\Phi_\nabla$.


If $Y$ is not connected, we define the reduced index at a point $\mu\in \Lambda$ as the sum of reduced indices over all connected components of $Y$.

More generally, let $H$ be a compact connected group acting on $Y$ and such that
$[H,H]$ acts trivially on $Y$. Let $\Scal_Y$ be an equivariant spin$^c$ structure on $Y$ with determinant bundle $\Lbb$.
For any $\mu\in \hgot^*$ such that $\mu([\hgot,\hgot])=0$, and  admissible for $H$, it is then possible to   define $\QS(Y_\mu)$.
 Indeed   eventually passing to a double cover of the torus
$H/[H,H]$ and translating by the square root of the action of $H/[H,H]$ on the fiber of $\Lbb$, we are reduced to the preceding case of the action of the torus $H/[H,H]$, and a $H/[H,H]$-equivariant spin$^c$ structure on $Y$.

Consider now the action of a  connected compact group $K$ on $M$.
Let $\sigma$ be a (relative interior) of a face of $\tgot^*_{\geq 0}$  which satisfies the following conditions
\begin{equation}\label{eq:C1-C2}
([\kgot_M,\kgot_M])=([\kgot_\sigma,\kgot_\sigma]),\hspace{0.5cm}\Phi_\nabla^{-1}(\sigma)\neq \emptyset.
\end{equation}


Let us explain how to compute the ``reduced'' index map
$\mu \to \QS(M_\mu)$ on the set  $\sigma\cap \{\Lambda+ \rho-\rho_\sigma\}$ that parameterizes the admissible orbits intersecting $\sigma$.
We work with the ``slice''  $Y$ defined by $\sigma$. The set $U_\sigma:=K_\sigma(\cup_{\sigma\subset \overline{\tau}}\tau)$ is an open
neighborhood of $\sigma$ in $\kgot^*_\sigma$ such that the open subset $KU_\sigma\subset \kgot^*$ is isomorphic to $K\times_{K_\sigma} U_\sigma$.
We consider the $K_\sigma$-invariant subset $Y=\Phi_\nabla^{-1}(U_\sigma)$.
The following lemma allows us to reduce the problem to the abelian case.
\begin{lemma}
$\bullet$ $Y$ is a non-empty submanifold of $M$ such that $KY$ is an open
susbset of $M$ isomorphic to $K\times_{K_\sigma}Y$.

$\bullet$ The Clifford module $\Scal_M$ on $M$ determines a Clifford module
$\Scal_{Y}$ on $Y$ with determinant line bundle $\Lbb_{Y}=\Lbb_M\vert_{Y}\otimes \C_{-2(\rho-\rho_\sigma)}$.
The corresponding moment map is $\Phi_\nabla\vert_{Y} - \rho+ \rho_\sigma$.

$\bullet$  The group $[K_\sigma,K_\sigma]$ acts trivially on $Y$ and on the bundle of spinors
$\Scal_{Y}$.
\end{lemma}

We thus consider $Y$ with action of $K_\sigma$, and Clifford bundle
$\Scal_{Y}$.
 If $\mu\in \sigma$ is admissible for $K$, then $\mu-\rho+\rho_\sigma\in \Lambda$ is admissible for $K_\sigma$.
The reduced space $M_\mu=\Phi_\nabla^{-1}(\mu)/K_\sigma$ is equal
to the reduced space $Y_{\mu-\rho+\rho_\sigma}$.
As $[K_\sigma,K_\sigma]$ acts trivially on $(Y,\Scal_Y)$, we are in the abelian case, and we define
$\QS(M_\mu):=\QS(Y_{\mu-\rho+\rho_\sigma})$.

\section*{Acknowledgments}
We wish to thank
  the Research in Pairs program at
Mathematisches Forschungsinstitut Oberwolfach (January 2014), which gave  us the opportunity to work on these questions.


\begin{thebibliography}{00}






\bibitem{CdS-K-T}{\sc A. Cannas da Silva}, {\sc Y. Karshon} and
{\sc S. Tolman}, Quantization of presymplectic manifolds and circle
actions, {\em Trans. Amer. Math. Soc.}  {\bf 352} (2000),  525-552.


\bibitem{du} M. Duflo, {\it Construction de repr\'esentations unitaires d'un groupe de Lie}, CIME, Cortona (1980).



\bibitem{Guillemin-Sternberg82} {\sc V. Guillemin} and {\sc S. Sternberg},
{\em Geometric quantization and multiplicities of group representations}, Invent. Math. {\bf 67} (1982), 515--538.

\bibitem{Grossberg-Karshon94} {\sc M. Grossberg} and  {\sc Y. Karshon}, {\em Bott towers, complete integrability, and the extended character of representations}, Duke Mathematical Journal {\bf 76} (1994), 23-58.

\bibitem{Grossberg-Karshon98} {\sc M. Grossberg} and  {\sc Y. Karshon}, {\em Equivariant index and the moment map for completely integrable torus actions},  Advances in Mathematics {\bf 133} (1998),  185-223.

\bibitem{Karshon-Tolman93}  {\sc Y. Karshon} and {\sc S. Tolman},
{\em The moment map and line bundles over presymplectic toric manifolds},
J. Differential Geom {\bf 38} (1993), 465-484.


\bibitem{Meinrenken98} {\sc E. Meinrenken}, {\em Symplectic surgery and the Spin\textsuperscript{c}-Dirac operator},
Advances in Math. {\bf 134} (1998),  240-277.

\bibitem{meinrenken-sjamaar} {\sc E. Meinrenken} and {\sc R. Sjamaar}, {\it Singular reduction and quantization}, Topology {\bf 38} (1999), 699-763.


\bibitem{pep-RR} {\sc P.-E. Paradan}, {\em Localization of the Riemann-Roch
character},   J. Functional Analysis  {\bf 187} (2001),  442--509.


\bibitem{pep-spin} {\sc P.-E.  Paradan}, {\it Spin-quantization commutes with reduction}, J. Symplectic Geometry {\bf 10} (2012), 389-422.

\bibitem{Tian-Zhang98} {\sc Y. Tian} and {\sc W. Zhang}, {\em An analytic proof of the geometric quantization conjecture of Guillemin-Sternberg},  Invent. Math. {\bf 132} (1998),  229--259.








\end{thebibliography}
\end{document}